# ON THE COHOMOLOGY OF $SL_2(Z[1/p])$

Alejandro Adem[†]
Mathematics Department
University of Wisconsin
Madison, WI 53706 USA

Nadim Naffah
Department of Mathematics
ETH–Zürich
Zürich CH–8092

## §0. INTRODUCTION

In this note we compute the integral cohomology of the discrete groups $SL_2(Z[1/p])$, $p$ a prime. According to Serre [S1] these are groups of virtual cohomological dimension 2. The method we use is to exploit the fact that these groups can be expressed as an amalgamation of two copies of $SL_2(Z)$ along the subgroup $\Gamma_0(p)$ of $2 \times 2$ matrices with lower left hand entry divisible by $p$. We first compute the cohomology of this virtually free group (using a tree on which it acts with finite isotropy and compact quotient), and then use the well–known Mayer–Vietoris sequence in cohomology to obtain our result. We assume that $p$ is an odd prime larger than 3. The cases $p = 2, 3$ must be treated separately; we discuss them at the end of the paper. We are grateful to the referee for his extremely useful remarks, and to J.-P. Serre for pointing out Proposition 3.1.

## §1. DOUBLE COSETS AND PERMUTATION MODULES

In this section we calculate certain double coset decompositions which will play a key rôle in our approach. Let $G = SL_2(F_p)$, and $B \subset G$ the subgroup consisting of all matrices with lower left hand entry equal to zero. It is easy to see that $B \cong Z/p \times_T Z/p - 1$, a semidirect product. We will denote by $C_2$, $C_4$ and $C_6$ the cyclic subgroups generated by the following three respective matrices of orders 2, 4 and 6:

$$a_2 = \begin{pmatrix} -1 & 0 \\ 0 & -1 \end{pmatrix}, \ a_4 = \begin{pmatrix} 0 & -1 \\ 1 & 0 \end{pmatrix}, \ a_6 = \begin{pmatrix} 0 & -1 \\ 1 & 1 \end{pmatrix}.$$

It is easy to check that the set of right cosets $B\backslash G$ decomposes as follows:

$$B\backslash G = Ba_4 \sqcup \Big( \bigsqcup_{x \in F_p} B \begin{pmatrix} 1 & 0 \\ x & 1 \end{pmatrix} \Big).$$

Now as $a_2$ is central, we obtain a double coset decomposition for G using $B$ and $C_2$ as follows:

$$G = Ba_4C_2 \sqcup \Big( \bigsqcup_{x \in F_p} B \begin{pmatrix} 1 & 0 \\ x & 1 \end{pmatrix} C_2 \Big).$$

[†] Partially supported by the NSF.

Applying the usual induction restriction formula we obtain

$$Z[G/C_2]\big|_B \cong (Z[B/C_2])^{p+1}. \tag{1.1}$$

Next we consider the double cosets using $C_4$. Note that we have

$$\begin{pmatrix} 1 & 0 \\ x & 1 \end{pmatrix} \begin{pmatrix} 0 & -1 \\ 1 & 0 \end{pmatrix} = \begin{pmatrix} 0 & -1 \\ 1 & -x \end{pmatrix}, \text{ and } \begin{pmatrix} 0 & -1 \\ 1 & -x \end{pmatrix} \begin{pmatrix} 1 & 0 \\ 1/x & 1 \end{pmatrix} = \begin{pmatrix} -1/x & -1 \\ 0 & -x \end{pmatrix}.$$

From this we conclude that if $x \neq -1/x$ or $0$, then

$$B \begin{pmatrix} 1 & 0 \\ x & 1 \end{pmatrix} \begin{pmatrix} 0 & -1 \\ 1 & 0 \end{pmatrix} = B \begin{pmatrix} 1 & 0 \\ -1/x & 1 \end{pmatrix}$$

and so

$$B \begin{pmatrix} 1 & 0 \\ x & 1 \end{pmatrix} C_4 = B \begin{pmatrix} 1 & 0 \\ x & 1 \end{pmatrix} \sqcup B \begin{pmatrix} 1 & 0 \\ -1/x & 1 \end{pmatrix}, \text{ and } Ba_4C_4 = B \sqcup Ba_4.$$

In each case the two cosets are permuted by the matrix $a_4$. In the case when $x^2 + 1 = 0$, then the coset is fixed under this action, and hence the associated coset is equal to the double coset. It is an elementary fact that the polynomial $t^2 + 1$ will have roots in $F_p$ (necessarily two distinct ones) if and only if $p \equiv 1 \mod(4)$. Using this and the induction restriction formula yields

$$Z[G/C_4]\big|_B \cong Z[B/s_1C_4s_1^{-1}] \oplus Z[B/s_2C_4s_2^{-1}] \oplus (Z[B/C_2])^{(p-1)/2} \quad \text{if } p \equiv 1 \mod(4)$$

and

$$Z[G/C_4]\big|_B \cong (Z[B/C_2])^{(p+1)/2} \quad \text{otherwise.} \tag{1.2}$$

The elements $s_1, s_2$ correspond to the two roots of the polynomial.

For $C_6$, we must look at the action of the matrix of order 3, $\begin{pmatrix} 0 & 1 \\ -1 & -1 \end{pmatrix}$ on the double cosets. In this case the orbit of the action will be of the form

$$B \begin{pmatrix} 1 & 0 \\ x & 1 \end{pmatrix}, \ B \begin{pmatrix} 1 & 0 \\ 1/(1-x) & 1 \end{pmatrix}, \ B \begin{pmatrix} 1 & 0 \\ (x-1)/x & 1 \end{pmatrix}$$

provided $x \neq 1$. A coset will be fixed if and only if $x^2 - x + 1 = 0$; given that $p > 3$, this will have roots in $F_p$ (necessarily two distinct ones) if and only if $p \equiv 1 \mod(3)$. As for $x = 1$, the corresponding coset gives rise to the singular orbit

$$B, \ B \begin{pmatrix} 1 & 0 \\ 1 & 1 \end{pmatrix}, \ Ba_4.$$

From this we can deduce the following decomposition (notation as before):

$$Z[G/C_6]\big|_B \cong Z[B/s_1 C_6 s_1^{-1}] \oplus Z[B/s_2 C_6 s_2^{-1}] \oplus \left(Z[B/C_2]\right)^{(p-1)/3} \quad \text{if } p \equiv 1 \bmod(3)$$

and

$$Z[G/C_6]\big|_B \cong \left(Z[B/C_2]\right)^{(p+1)/3} \quad \text{otherwise.} \tag{1.3}$$

## §2 THE COHOMOLOGY OF $\Gamma_0(p)$

The subgroup $\Gamma_0(p) \subset SL_2(Z)$ is defined by

$$\Gamma_0(p) = \{\begin{pmatrix} a & b \\ c & d \end{pmatrix} \in SL_2(Z) \mid c \equiv 0 \bmod(p)\}.$$

If $\Gamma(p)$ denotes the level $p$ congruence subgroup, then clearly $\Gamma_0(p)$ can be expressed as an extension

$$1 \to \Gamma(p) \to \Gamma_0(p) \to B \to 1.$$

Recall [S1] that $SL_2(Z)$ acts on a tree $T$ with finite isotropy, and quotient a single edge. As $\Gamma(p)$ is torsion free, it acts freely on this tree, and so $G = SL_2(F_p)$ acts on the finite graph $T/\Gamma(p)$. The isotropy subgroups of this action are precisely $C_4$ and $C_6$ for the vertices and $C_2$ for the edge. Let $EB$ denote the universal $B$–space, then clearly using the projection $\pi : \Gamma_0(p) \to B$, $\Gamma_0(p)$ can be made to act diagonally on $EB \times T$, with trivial isotropy. As this space is contractible, its quotient under this action has the same homotopy type as the classifying space $B\Gamma_0(p)$, and so we have $B\Gamma_0(p) \cong EB \times_B T/\Gamma(p)$.

Let $C^*$ denote the cellular cochains on the $B$–CW complex $T/\Gamma(p)$; then it is clear from the above that

$$C^0 \cong Z[G/C_4]\big|_B \oplus Z[G/C_6]\big|_B \quad \text{and that} \quad C^1 \cong Z[G/C_2]\big|_B.$$

In this situation, there is a spectral sequence (see [B]) converging to $H^*(\Gamma_0(p), Z)$, with $E_1^{p,q} \cong H^q(B, C^p)$, which degenerates into a long exact sequence:

$$\ldots \to H^i(\Gamma_0(p), Z) \to H^i(B, C^0) \to H^i(B, C^1) \to H^{i+1}(\Gamma_0(p), Z) \to \ldots \tag{2.1}$$

where the middle arrow is induced by the coboundary map $\delta$ on $C^*$.

As a first application of the long exact sequence (2.1) we obtain

**Proposition 2.2**

Under the above conditions, $H^1(\Gamma_0(p), Z) \cong (Z)^{N(p)}$ where

$$N(p) = \begin{cases} (p-7)/6, & \text{if } p \equiv 1 \bmod(12); \\ (p+1)/6, & \text{if } p \equiv 5 \bmod(12); \\ (p-1)/6, & \text{if } p \equiv 7 \bmod(12); \\ (p+7)/6, & \text{if } p \equiv 11 \bmod(12). \end{cases}$$

**Proof:** The sequence (2.1) starts as

$$0 \to Z \to (C^0)^B \to (C^1)^B \to H^1(\Gamma_0(p), Z) \to H^1(B, C^0).$$

Recall that $H^1(B, C^0) = 0$, as $C^0$ is a permutation module. Hence calculating ranks completes the proof. ∎

To compute the remaining cohomology groups we first switch to the associated projective group $P\Gamma_0(p)$. Note that there will be a situation analogous to that for the original group, except that throughout we must divide out by the central $Z/2$. Note that if $PB$ is the associated group for $B$, then $C^1$ will now be a free $PB$–module. Hence the corresponding long exact sequence degenerates to yield the isomorphism $H^{2i}(P\Gamma_0(p), Z) \cong H^2(PB, C^0)$ for all $i \geq 0$, and the fact that all its odd dimensional cohomology (except in dimension 1) is zero. This is summarized in

**Proposition 2.3**

For any integer $i \geq 1$, we have that

$$H^{2i}(P\Gamma_0(p), Z) \cong \begin{cases} Z/6 \oplus Z/6, & \text{if } p \equiv 1 \ mod(12); \\ Z/2 \oplus Z/2, & \text{if } p \equiv 5 \ mod(12); \\ Z/3 \oplus Z/3, & \text{if } p \equiv 7 \ mod(12); \\ 0, & \text{if } p \equiv 11 \ mod(12). \end{cases}$$

and $H^{2i+1}(P\Gamma_0(p), Z) = 0$. ∎

Note that $H^1(P\Gamma_0(p), Z) \cong H^1(\Gamma_0(p), Z)$.

Next we apply the spectral sequence over $Z$ associated to the central extension

$$1 \to C_2 \to \Gamma_0(p) \to P\Gamma_0(p) \to 1.$$

Note that the interesting cases are if $p \equiv 1, 5 \mod (12)$, and that the 3–torsion plays no role. As the group has periodic cohomology and is virtually free, it suffices to compute $H^2$ and $H^3$. In total degree 2 we simply have the contributions from $H^2(C_2, Z) \cong Z/2$ and $H^2(P\Gamma_0(p), Z)_{(2)} \cong Z/2 \oplus Z/2$. As 4–torsion must appear (there will be a subgroup of that order), we conclude that $H^2(\Gamma_0(p), Z)_{(2)} \cong Z/4 \oplus Z/2$.

In total degree 3, we only have one term: $H^1(P\Gamma_0(p), H^2(C_2, Z)) \cong (Z/2)^{N(p)+2}$. However, it is not hard to see that the map induced by the quotient in cohomology, $H^4(P\Gamma_0(p), Z) \to H^4(\Gamma_0(p), Z)$ must be zero. This can be proved by comparing the two long exact sequences described above and using the corresponding fact for the map induced by the quotient $Z/4 \to Z/2$. The only possible differential on this horizontal edge group is $d_3 : E_3^{1,2} \to E_3^{4,0}$, hence it must have an image of 2–primary rank 2. We conclude that $E_\infty^{1,2} \cong (Z/2)^{N(p)} \cong H^3(\Gamma_0(p), Z)$, and we obtain

**Theorem 2.4**

For any integer $i \geq 1$, we have

$$H^{2i}(\Gamma_0(p), Z) \cong \begin{cases} Z/12 \oplus Z/6, & \text{if } p \equiv 1 \bmod(12); \\ Z/4 \oplus Z/2, & \text{if } p \equiv 5 \bmod(12); \\ Z/3 \oplus Z/6, & \text{if } p \equiv 7 \bmod(12); \\ Z/2, & \text{if } p \equiv 11 \bmod(12) \end{cases}$$

and

$$H^{2i+1}(\Gamma_0(p), Z) \cong (Z/2)^{N(p)}.$$

## §3. CALCULATION OF THE COHOMOLOGY

To begin we recall that aside from the natural inclusion we also have an injection $\rho : \Gamma_0(p) \to SL_2(Z)$, given by

$$\begin{pmatrix} a & b \\ c & d \end{pmatrix} \mapsto \begin{pmatrix} a & pb \\ p^{-1}c & d \end{pmatrix}.$$

Using these two imbeddings, we can construct the amalgamated product [S1]

$$SL_2(Z[1/p]) \cong SL_2(Z) *_{\Gamma_0(p)} SL_2(Z).$$

In addition we have that

$$\overline{H}^r(SL_2(Z), Z) \cong \begin{cases} Z/12 & \text{if r is even;} \\ 0 & \text{if r is odd.} \end{cases}$$

We can identify $\rho^*$ with the ordinary restriction map.

Using the Mayer–Vietoris sequence associated to an amalgamated product we see that $H^1(SL_2(Z[1/p], Z) = 0$ and that we have exact sequences

$$0 \to Z^{N(p)} \to H^2(SL_2(Z[1/p]), Z) \to Z/12 \oplus Z/12 \to H^2(\Gamma_0(p), Z) \to H^3(SL_2(Z[1/p]), Z) \to 0$$

and

$$0 \to (Z/2)^{N(p)} \to H^{2i}(SL_2(Z[1/p], Z) \to Z/12 \oplus Z/12 \to H^{2i}(\Gamma_0(p), Z) \to H^{2i+1}(SL_2(Z[1/p]), Z) \to 0.$$

The cohomology will evidently be 2-fold periodic above dimension 2, which is in fact the virtual cohomological dimension of $SL_2(Z[1/p])$.

We will need the following result, which is due to J.-P. Serre [S2]. It can also be proved using an explicit presentation for the group, described in [BM].

**Proposition 3.1:**

$$H_1(SL_2(Z[1/p], Z) \cong \begin{cases} Z/3 & \text{if p=2;} \\ Z/4 & \text{if p=3;} \\ Z/12 & \text{otherwise.} \end{cases}$$

Hence we have that for $p > 3$, $H^2(SL_2(Z[1/p], Z) \cong (Z)^{N(p)} \oplus Z/12$.

Let $A(p)$ denote the number $12/|Q(p)|$, where $Q(p)$ is the largest cyclic subgroup in $H^2(\Gamma_0(p), Z)$. Then, from the fact that the restriction from the cohomology of $SL_2(Z)$ to that of its cyclic subgroups factors through $\Gamma_0(p)$, we deduce that the sequence above simplifies to yield

(3.2) $$0 \to (Z/2)^{N(p)} \to H^{2i}(SL_2(Z[1/p], Z) \to Z/12 \oplus Z/A(p) \to 0$$

and

$$H^{2i+1}(SL_2(Z[1/p]), Z) \cong H^{2i}(\Gamma_0(p), Z)/Q(p).$$

Moreover, from our previous calculation for $\Gamma_0(p)$, we have that

$$Q(p) = \begin{cases} Z/12 & \text{if } p \equiv 1 \bmod(12); \\ Z/4 & \text{if } p \equiv 5 \bmod(12); \\ Z/6 & \text{if } p \equiv 7 \bmod(12); \\ Z/2 & \text{if } p \equiv 11 \bmod(12). \end{cases}$$

It remains only to determine precisely what this extension (3.2) looks like. We need only be concerned with the 2–primary component. Recall that $H^{2i}(SL_2(Z[1/p]), Z)$ is a quotient of $H^2(SL_2(Z[1/p]), Z)$, as $SL_2(Z[1/p])$ is a group of virtual cohomological dimension 2, which has 2-fold periodic cohomology. This can also be explained by saying that the sequence in high even dimensions can be identified with the corresponding sequence in 2-dimensional Farrell cohomology (see [B]). Mapping one sequence into the other, we see that a $Z/4$ summand must split off for all values of $p > 3$. This means that the sequence *will split* for $p \equiv 1, 5 \bmod(12)$. For the remaining cases $p \equiv 7, 11 \bmod(12)$, it remains to solve the extension problem after splitting off the $Z/4$ summand. However, from our knowledge of $H^2$, we know that the dimension of $H^{2i} \otimes Z/2$ can be at most $N(p) + 1$. We infer that the reduced extension *does not split*, and must necessarily have a $Z/4$ summand present. Hence we have the following complete calculation:

**Theorem 3.3**

Let $p$ be an odd prime larger than 3, then $H^1(SL_2(Z[1/p]), Z) \cong 0$, and

$$H^2(SL_2(Z[1/p]), Z) \cong \begin{cases} Z^{(p-7)/6} \oplus Z/12 & \text{if } p \equiv 1 \bmod(12); \\ Z^{(p+1)/6} \oplus Z/12 & \text{if } p \equiv 5 \bmod(12); \\ Z^{(p-1)/6} \oplus Z/12 & \text{if } p \equiv 7 \bmod(12); \\ Z^{(p+7)/6} \oplus Z/12 & \text{if } p \equiv 11 \bmod (12). \end{cases}$$

For $i \geq 2$, we have

$$H^{2i}(SL_2(Z[1/p]), Z) \cong \begin{cases} (Z/2)^{(p-7)/6} \oplus Z/12 & \text{if } p \equiv 1 \bmod(12); \\ (Z/2)^{(p+1)/6} \oplus Z/12 \oplus Z/3 & \text{if } p \equiv 5 \bmod(12); \\ (Z/2)^{(p-7)/6} \oplus Z/12 \oplus Z/4 & \text{if } p \equiv 7 \bmod(12); \\ (Z/2)^{(p+1)/6} \oplus Z/12 \oplus Z/12 & \text{if } p \equiv 11 \bmod(12), \end{cases}$$

and

$$H^{2i-1}(SL_2(Z[1/p]), Z) \cong \begin{cases} Z/6 & \text{if } p \equiv 1 \ mod(12); \\ Z/2 & \text{if } p \equiv 5 \ mod(12); \\ Z/3 & \text{if } p \equiv 7 \ mod(12); \\ 0 & \text{if } p \equiv 11 \ mod(12). \end{cases}$$

∎

Of the two remaining cases ($p = 2, 3$) the second one can be done in a manner totally analogous to what we have presented. Specifically we have that

$$H^i(\Gamma_0(3), Z) \cong \begin{cases} Z & \text{if } i = 0, 1; \\ Z/6 & \text{if } i \text{ is even}; \\ Z/2 & \text{if } i > 1 \text{ is odd}. \end{cases}$$

We obtain that $H^2(SL_2(Z[1/3], Z) \cong Z \oplus Z/4$, $H^{2i+1}(SL_2(Z[1/3]), Z) \cong 0$, and it only remains to deal with the extension

$$0 \to Z/2 \to H^{2i}(SL_2(Z[1/p], Z) \to Z/12 \oplus Z/2 \to 0.$$

As before the $Z/12$ must split off, and by rank considerations an extra $Z/4$ summand must appear. To summarize, we have

$$H^i(SL_2(Z[1/3]), Z) \cong \begin{cases} 0 & \text{if } i \text{ is odd}; \\ Z \oplus Z/4 & \text{if } i = 2; \\ Z/12 \oplus Z/4 & \text{if } i = 2j, \ j > 1. \end{cases}$$

The case $p = 2$ is complicated by the fact that $\Gamma(2)$ is not torsion–free. However, one can still make use of the associated projective groups; $P\Gamma(2)$ is free of rank 2, and there is an extension

$$1 \to P\Gamma(2) \to P\Gamma_0(2) \to Z/2 \to 1.$$

Analyzing the action on the corresponding graph, it is not hard to show that

$$H^*(P\Gamma_0(2), Z) \cong \begin{cases} Z & \text{if } i \text{ is 0 or 1}; \\ Z/2 & \text{if } i \text{ is even}; \\ 0 & \text{otherwise}. \end{cases}$$

Then, using the central extension, it is direct to show that

$$H^i(\Gamma_0(2), Z) \cong \begin{cases} Z & \text{if } i \text{ is 0 or 1}; \\ Z/4 & \text{if } i \text{ is even}; \\ Z/2 & \text{if } i \text{ is odd}, i > 1. \end{cases}$$

Using the Mayer–Vietoris sequence as before, we obtain that $H^2(SL_2(Z[1/2]), Z) \cong Z \oplus Z/3$, $H^i(SL_2(Z[1/2]), Z) = 0$ if $i$ is odd, and a short exact sequence

$$0 \to Z/2 \to H^{2i}(SL_2(Z[1/2]), Z)) \to Z/12 \oplus Z/3 \to 0.$$

In this case, we do not know that the $Z/12$ summand must split off. Looking at the 2-primary part, we see that it can have *at most* one cyclic summand. The only possibility is $Z/8$, and we have

$$H^i(SL_2(Z[1/2]), Z) \cong \begin{cases} 0 & \text{if } i \text{ is odd}; \\ Z \oplus Z/3 & \text{if } i = 2; \\ Z/24 \oplus Z/3 & \text{if } i = 2j, \ j > 1. \end{cases}$$

*Remarks*:

Note that this last group has no finite subgroups of order eight, which makes its cohomology rather interesting. We are grateful to Hans–Werner Henn for pointing out the correct cohomology of this group. Also we would like to point out that Naffah [N] has calculated the 3-adic component of the Farrell cohomology $\widehat{H}^*(SL_2(Z[1/N]), Z)$ for *any* integer $N$. This of course can be used to recover our calculations at $p = 3$ in dimensions larger than 2. Moss [M] has computed the rational cohomology of $SL_2(Z[1/N])$ which again can be used to recover part of our results. The general calculation of $H^*(SL_2(Z[1/N]), Z)$ seems to be a rather complicated but interesting open problem. We refer to [H] for more on this.